\documentstyle[amssymb,12pt]{article}
\font\tmsb=msbm10 at12pt
\font\smsb=msbm7
\font\ssmsb=msbm5
\newfam\msbfam
\textfont\msbfam=\tmsb
\scriptfont\msbfam=\smsb
\scriptscriptfont\msbfam=\ssmsb
\def \mth {\fam\msbfam}
\def \Mth#1 {{\mth #1}}

\def \NN {{\Bbb {N}}}

\newcommand \bd {\begin {displaymath}}
\newcommand \ed {\end {displaymath}}

\newcommand{\Si}{\Sigma}

\newtheorem {defi} {Definition} [section]
\newtheorem {lem} {Lemma} [section]
\newtheorem {th} {Theorem} [section]

\newtheorem {form} {Formula} [section]

\title{The moduli space ${\cal M}_n(\Si)$ of stable fiber bundles over a 
 compact Riemann surface}
\author{A.Balan\\
Ecole Polytechnique\\
Centre de mathematiques\\
 UMR 7640 of CNRS\\
F-91128 PALAISEAU Cedex\\
{\it email}: balan@math.polytechnique.fr}
\date{}
\begin{document}
\maketitle

\abstract{ The moduli space of holomorphic 
fiber bundles  ${\cal M}_n(\Si)$ over a compact 
Riemann surface $\Si$ is considered.
 A formula for the regularised determinant
 and an other for the symplectic form at trivial bundle are proposed.}

\newpage

\section{Introduction}
 The moduli space ${\cal M}_n(\Si)$ of holomorphic stable fiber bundles
 of rank $n$ and trivial determinant over a compact Riemann surface
 and of genus $g \geq 2$ possesses usually three descriptions.
 The first one is given by the {\it symplectic quotient} \cite{GS}
 of the space of the unitary connections ${\cal A}_{\Si}$ over the trivial
 fiber bundle when the moment action given by the curvature is considered.
 The second one is a double algebraic quotient called {\it field} \cite{So}
 which represents the holomorphic fiber bundles.
 The third one is furnished by a {\it quotient of the set of irreductible
 homomorphisms} \cite{NS} of the Poincar\'e group of the Riemann surface
$\pi_1 (\Si)$ in the group of unitary symmetries of rank $n$, over the
 action of the same group. This space is a non-compact complex variety.
 The regularised determinant gives the 
 corresponding metric, called Quillen metric,
 a formula for the regularised determinant is proposed
 in the formula \ref{detreg}.
 It exits, morover, a symplectic form over this space which is a kaehlerian
 variety \cite{Hi}. Here is presented also a formula for the symplectic form
 at trivial fiber bundle in theorem \ref{symp} by mean of a reduction of
 the form by the harmonic ones over the Riemann surface.
\section{The cocycles and the operators $\overline \partial$}
\subsection{The connection with a $1$-cocycle, $f(z)$}
 Let a compact Riemann surface $\Si$ be, with a point $p$ and a
 holomorphic injection of the unitary disc of the complex plan of the surface:
   $$ i : \ \ (D,0) \hookrightarrow ( \Si , p).$$
 Let $\varrho$ a partition of the unity be defined for the disc
 and the surface minus a point.
 Let a holomorphic $1$-cocycle be with valus in the Lie algebra,
 $Lie(Sl_n)({\mathbb C})$ and defined over the disc minus a point:
 $$ f(z) : [ D -\{ 0 \}] \rightarrow Lie(Sl_n)({\mathbb C}) .$$
 The holomorphic cocycle can then be considered with values in the Lie group,
 $Sl_n({\mathbb C})$ \cite{DG}:
 $$ exp(f)(z): [D-\{ 0 \}] \rightarrow Sl_n({\mathbb C}).$$
 Now, a differential operator is defined associated with these datas,
 acting over the trivial fiber bundle of rank $n$ and the solutions
 define a holomoprhic fiber bundle equivalent with the fiber bundle which
 is obtained with a holomorphic cocycle. The following operator is considered:
 $$ \overline {\partial} + f(z) \overline {\partial} \varrho.$$
 This operator is well defined, the partition of unity being constant
 near the point $p$ and outside the disc $D$. The connection
 is defined taking the adjoint for the trivial metric:
  $$    d + f(z) \overline {\partial} \varrho  - f(z)^*{\partial} \varrho. $$
\subsection{The isomorphism ${\cal C}^{\infty}$ with the trivial fiber
 bundle}
 The holomorphic fiber bundle is associated with a $1$-cocycle:
\medskip
$$ ([\Si-D]. {\mathbb C}^n) \coprod (D. {\mathbb C}^n) / R,$$
where $R$ is the equivalence relation given by the holomorphic cocycle:
 $$ \forall z \in S^1; \forall s \in {\mathbb C}^n;$$
 $$ \ (z,s) R (z,[\exp(f)(z)]s) .$$
 There is an injection of a section ${\cal C}^{\infty}$ of the fiber bundle
 in a section of the trivial fiber bundle by mean of the following
 applications, which define the injection:
 $$ f_1 : \ (z,s) \in ([\Si-D].  {\mathbb C}^n)
 \mapsto (z,( \exp( -\varrho f(z)))s).$$
 $$ f_2 :\ (z,s) \in (D.{\mathbb C}^n)
 \mapsto (z,(\exp(-(1-\varrho)f(z)))s).$$ 
$
\begin{array}{crcrcr}
 & ( [\Si-D]. {\mathbb C}^n) \coprod (D. {\mathbb C}^n) / R & \rightarrow & 
 (\Si-D). {\mathbb C}^n\\
   &  \downarrow                        &          & \downarrow^{f_1} \\
      &   D. {\mathbb C}^n   & \stackrel {f_2} {\rightarrow}  &  \   \Si.
 {\mathbb C}^n   \\
\end{array}
$
$$([\Si-D]. {\mathbb C}^n) \coprod (D. {\mathbb C}^n) / R
 \rightarrow  \Si. {\mathbb C}^n,$$
by the universal property of the product in the category of the
 fiber bundles ${\cal C}^{\infty}$ over the Riemann surface $\Si$.
 This arrow is a trivialisation ${\cal C}^{\infty}$ of the fiber bundle.
 The holomorphic local sections of the fiber bundle defined with the cocycle
 are taken by the above defined application in these which verify over
 ${\Si}^*$:
$$ [ \exp (- \varrho f(z) ) \circ \overline {\partial} \circ
\exp ( \varrho f(z) ) ] \ s  =
 \overline {\partial} s + f(z) \overline {\partial} \varrho \  s = 0.$$
 So, a connection is obtained canonicaly
 associated with the $1$-cocycle, taking the adjoint by the metric.
\subsection{The distributions}
 The differential operator which is obtained is dependant of the choice
 of a partition of the unity.
 To avoid the choice, the partition of unity must go toward
 the caracteristic function of the disc, in the sens of the distributions
 \cite{Sc}, so an operator of the type is obtained:
 $$ \overline \partial + [ f(z) \delta_{S^1}  / 2 \pi] d \overline z : $$
 $$   {\Gamma }^0  (\Si, {\mathbb C}^n) \longrightarrow 
\Gamma^{0,1} (\Si, {\mathbb C}^n)  \hat \otimes_{{\cal C}^{\infty}(\Si)}
 {\cal D}'(\Si, {\mathbb C}).$$
 Where ${\cal D}'(\Si, {\mathbb C})$ is the space of distributions
 and  $\delta_{S^1}/ 2 \pi$, is a distribution associated with the circle.
 The tensor product is taken with the ring of functions
 ${\cal C}^{\infty}$ over $\Si$. So the operator of derivation plus
 an operator with values in the distributions over the Riemann surface $\Si$ is
 obtained.
\bigskip
 The whole connection is obtained taking the adjoint of the operator
 for the metric over the trivial fiber bundle;
 it defines the part of type $ (1,0)$ with compatibility
 with the condition of unitarity of the connection:
$$ d +  f(z) \overline \partial \varrho - f(z)^* \partial \varrho.$$
 To avoid the choice of a partition of unity, it must
 go towards the caracteristic function of the disc;
 then an operator with values in the distributions of the following type
 is obtained:
$$d + [ f(z) \delta_{S^1}  d \overline z - f(z)^* \delta_{S^1} dz]/ 2 \pi :$$
$$   {\Gamma}^{0} (\Si, {\mathbb C}^n) \longrightarrow 
\Gamma^{1} (\Si, {\mathbb C}^n) \hat \otimes_{{\cal C}^{\infty}(\Si)} {\cal D}'(\Si,
 {\mathbb C}) .$$
\section{The tangent space, $T {\cal M}_n(\Si)$}
 With the different descriptions of the space of moduli of the
 fiber bundles of the Riemann surface $\Si$, the tangent space of the
 moduli, $T {\cal M}_n(\Si)$, can be presented.

\medskip
\noindent
 In the first case, some complex forms of type $(0,1)$ over the Riemann
 surface $\Si$, taken as elements of the tangent space to the connections
 space, ${\cal A}_{\Si}$.

\medskip
\noindent
 In the second case, the Lie algebra of the groups of infinite dimension
 which is considered $ Lie( Sl_n)((z))$, and the injections.

\medskip
\noindent
 In the third case:
 $$Hom( \pi_1 (\Si), Lie(SU_n)( {\mathbb C})).$$
\section{The Quillen metric}
 The Quillen metric is a metric over the determinant fiber bundle
 of the moduli space, and is given by regularised determinants.
 Let a Riemann surface $\Si$ be, with boundary $ S^1 $.
 The derivation of the regularised determinant is calculated
\cite{Al}for a family of connections for a trivial holomorphic structure
 and a boundary condition, for example of Dirichlet type.
\subsection{The regularised determinants}
 The regularised determinant of a infinite positiv autoadjoint operator $P$
 over a hibert space of infinite dimension is defined by the proper values
 $\lambda_i$, $i \in  {\mathbb N}^*$. 
 First the theta function of the operator is defined:
$$ \theta_P (t) = \sum_{i \in {\mathbb N}^*}
 \exp(- t \lambda_i) \ ,\  t \in {{\mathbb R}_+}^*.$$ 
 Where the zero proper values are taken off.
 This series of positiv terms is supposed convergent for all $t$.
 The theta function is supposed to admit an asymptotic development
 when $t$ goes towards zero in the following form:
$$ \theta_P (t) = \sum_{n} a_n t^n .$$
 Where $n$ takes its values in the arithmetic
 increasing sequence of rationals.
The zeta function of the operator $\zeta_P$, which is the Mellin
 transform of the theta function, is then well defined and
 has an analytic prologation over the complex plan with poles in the
 numbers $-n+1$.
 In the case where zero is a regular value, the regularised determinant
 is defined as:
 $$ det'(P)= \exp(- \zeta'(0)). $$
 The definition can be applied in the case of a Laplacian
 with a connection over the above fiber bundle of a Riemann surface.
 The theta function is well defined from the fact that the heat kernel
 gives a trace operator.
 The asymptotic development is the one of Minakshisundaram-Pleijel.
\subsection{The case of rank $1$}
 First the case of a one dimensional fiber bundle is considered.
 The family of connections over the trivial fiber bundle is of the form:
 $$ A_t^{0,1} = \bar \partial + t \bar \partial f 
= \exp(-tf) \circ \bar \partial \circ \exp(tf) .$$
 Where $f$ is a real function over the Riemann surface.
 The part of type $(1,0)$ being given by the adjoint
 by mean of the canonical metric over the trivial fiber bundle.
 The considered Laplacian are then:
 $$ \Delta_t = \exp(tf) \circ \bar {\partial}^* \circ \exp(-2tf)
 \circ \bar \partial \circ \exp(tf) .$$ 
 A Dirichlet condition is imposed over the boundary to obtain 
 a positiv autoadjoint operator.
 The regularised determinant of the family of operators is then 
 derivable. An expression for the derivation is seeked.
 The following formula is taken:
 $$ \int_{\epsilon} ^{\infty}  \theta (t) \frac{dt}{t} = \int_0^{\infty}
 \theta (t) \frac{dt}{t}
+ \sum_{n \leq 0}  ( \frac{a_n}{n} - \frac{ a_n \epsilon^n}{n})
 + a_0  \ log(\epsilon)+ \int_{\epsilon} ^1 r(t) \frac{dt}{t} . $$
 Where $r$ is the theta function  without divergent terms in zero.
 The following equality is deduced:
 $$ \int_{\epsilon}^{\infty} \theta (t) \frac{dt}{t} = \zeta ' (0) 
 + \sum_{n \leq 0}  a_n \frac{\epsilon^n}{n} + a_0 log(\epsilon) +
 \int_{0}^{\epsilon} r(t) \frac{dt}{t}. $$
 Then, a derivation with respect to the the parameter of the family of
 operators, $s$:
 $$ \int_{\epsilon}^{\infty}\frac{ d}{ds} \theta (t) \frac{dt}{t} =
 \frac{d}{ds} \zeta ' (0) + \sum_{n \leq 0} \frac{d}{ds} a_n
 \frac{\epsilon^n}{n} + \frac{d}{ds} a_0 log(\epsilon) +
 \int_{0}^{\epsilon} \frac{d}{ds} r(t) \frac{dt}{t}. $$
 So:
$$\frac{ d}{ds} \zeta '(0) = PF \int_{\epsilon}^{\infty}
 \frac{d}{ds} \theta (t) \frac{dt}{t}.$$
 Where $PF$ is the finite part when $\epsilon$ near zero.
 $ \theta (t) = tr( \exp(-t \Delta )) $, 
 where the heat kernel associated with the Laplacian is considered,
 it is finaly obtained:
 $$ \frac{d}{ds} \zeta '(0) = 
 PF tr( (\frac{d}{ds} \Delta ) \Delta^{-1} \exp( - \epsilon \Delta)).$$
 Let now $| \Phi_i \rangle, i \in {\mathbb N}^*$,
 an orthogonal basis of proper vectors of the positiv autoadjoint operator
 $\Delta$ for the proper values $\lambda_i$; for this basis, the trace is
 expressed as:
 $$ tr( (\frac{d}{ds} \Delta) \Delta^{-1} 
\exp( - \epsilon \Delta) = \sum_{i}
 \langle \Phi_i| (\frac{d}{ds} \Delta ) \Delta^{-1} 
\exp( - \epsilon \Delta) | \Phi_i \rangle = $$
 $$
= \sum_{i}  \langle \Phi_i| (\frac{d}{ds} \Delta )
| \Phi_i \rangle \lambda_i^{-1} \exp( - \epsilon \lambda_i ).$$
 As the $|\Phi_i\rangle$ are proper vectors of the operator $\Delta$, it holds:
$$  \Delta |\Phi_i\rangle = \lambda_i|\Phi_i\rangle .$$
 Derivating, it is obtained that:
$$ (\displaystyle {\frac{d}{ds} } \Delta )|\Phi_i\rangle + 
\Delta (\frac{d}{ds}|\Phi_i\rangle) =
(\frac{d}{ds} \lambda_i)|\Phi_i\rangle  +\lambda_i
 (\frac{d}{ds}|\Phi_i\rangle).$$
 The variations of the operators $\Delta$ are then in the form:
$$ \frac{d}{ds} \Delta = f \exp(sf)  \circ {\bar \partial}^* \circ \exp(-2sf)
 \circ \partial \circ \exp(sf) +  $$ 
$$ + \exp(sf)  \circ {\bar \partial}^* \circ (-2f) \exp(-2sf) \circ
 \bar \partial \circ \exp(sf) + $$ 
$$ + \exp(sf)  \circ {\bar \partial}^* \circ \exp(-2sf) \circ
 \bar \partial \circ \exp(sf) f .$$
 This expression is put in the calculation of the trace and using the fact
 that the vectors $|{\Phi}_i\rangle$ are unitary, it is obtained:
$$ 
 PF tr( (\frac{d}{ds} \Delta ) \Delta^{-1} \exp( - \epsilon \Delta)).$$
 Let now $| \Phi_i \rangle, i \in {\mathbb N}^*$,
 an orthogonal basis of proper vectors of the positiv autoadjoint operator
 $\Delta$ for the proper values $\lambda_i$; for this basis, the trace is
 expressed as:
 $$ tr( (\frac{d}{ds} \Delta) \Delta^{-1} 
\exp( - \epsilon \Delta) = \sum_{i}
 \langle \Phi_i| (\frac{d}{ds} \Delta ) \Delta^{-1} 
\exp( - \epsilon \Delta) | \Phi_i \rangle = $$
 $$
= \sum_{i}  \langle \Phi_i| (\frac{d}{ds} \Delta )
| \Phi_i \rangle \lambda_i^{-1} \exp( - \epsilon \lambda_i ).$$
 As the $|\Phi_i\rangle$ are proper vectors of the operator $\Delta$, it holds:
$$  \Delta |\Phi_i\rangle = \lambda_i|\Phi_i\rangle .$$
 Derivating, it is obtained that:
$$ (\displaystyle {\frac{d}{ds} } \Delta )|\Phi_i\rangle + 
\Delta (\frac{d}{ds}|\Phi_i\rangle) =
(\frac{d}{ds} \lambda_i)|\Phi_i\rangle  +\lambda_i
 (\frac{d}{ds}|\Phi_i\rangle).$$
 The variations of the operators $\Delta$ are then in the form:
$$ \frac{d}{ds} \Delta = f \exp(sf)  \circ {\bar \partial}^* \circ \exp(-2sf)
 \circ \partial \circ \exp(sf) +  $$ 
$$ + \exp(sf)  \circ {\bar \partial}^* \circ (-2f) \exp(-2sf) \circ
 \bar \partial \circ \exp(sf) + $$ 
$$ + \exp(sf)  \circ {\bar \partial}^* \circ \exp(-2sf) \circ
 \bar \partial \circ \exp(sf) f .$$
 This expression is put in the calculation of the trace and using the fact
 that the vectors $|{\Phi}_i\rangle$ are unitary, it is obtained:
$$ 
 \sum_{i}  \langle \Phi_i| (\frac{d}{ds} \Delta )
| \Phi_i \rangle \lambda_i^{-1} \exp( - \epsilon \lambda_i) = $$
$$
= 2  \sum_{i}  \langle \Phi_i| f | \Phi_i \rangle 
\exp( - \epsilon \lambda_i) -  2 \sum_{i} \langle \Psi_i| f |\Psi_i \rangle 
\exp( - \epsilon \lambda_i).$$
 Where $|\Psi_i\rangle = \exp(-sf) \circ \partial \circ
 \exp(sf) |\Phi_i\rangle / \sqrt{\lambda_i}$
 are the orthogonal vectors and the propeer vectors for the proper
 values $\lambda_i$ of the operator: 
$$ \Delta^- = \exp(-sf)  \circ {\bar \partial}
 \circ \exp(2sf) \circ {\bar \partial}^* \circ \exp(-sf),$$
 acting over the $1$-forms over the surface.
 The variations of the determinant are then:
{\form .
\begin{equation}
 \frac{ d}{ds} \zeta '(0)= PF \  2tr(f \exp( - \epsilon \Delta) )
 -  PF \ 2 tr(f \exp(- \epsilon \Delta^-)). \label{detreg}
\end{equation}
}
 To obtain an explicite expression of the derivation of the
 regularised determinant, the first coefficients of the asymptotic
 development of the heat kernel must be calculated when the parameter
 goes toward zero by mean of the Seeley methods with boundary terms.
\subsection{The case of higher ranks}
 Let the trivial fiber bundle be of rank $n$ and a connection
 put in the following form:
 $$ A^{0,1} = \bar \partial + \alpha^{-1} \bar \partial {\alpha}. $$ 
Where $ \alpha $ is an inversible matrix.
 A polar decomposition is applied for the matrix:
 $$ \alpha =  H \circ U .$$ 
Where $U$ is unitary and the $H$ hermitian.
The connection is then put in the following form:
$$ U^* \circ  H^{-1} \circ \bar \partial \circ H \circ U .$$
And the associated Laplacian is:
$$ U^* \circ H \circ {\bar \partial}^* \circ H^{-2 } 
\circ \bar \partial \circ H \circ U. $$
 For the calculation of the regularised determiant, a positiv hermitian
 matrix is taken and the parametrised family of Laplacian
 is considered with the real $t$:
$$ H^t \circ {\bar \partial}^* \circ H^{-2t } \circ \bar \partial \circ H^t .$$
 So a similary calculation is done.

\medskip
\noindent
 An expression for the derivation of the regularised determinant
 has been obtained in the case of a conjugated connection over a Riemann
 surface with boundary.
 The derivation could be integrated \cite{DP}
 for an exatc formula of the regularised determinant.

\section{The symplectic form of the moduli space ${\cal M}_n(\Si)$}
\subsection{The symplectic form of the space of moduli, $\omega$ }
 When the holomorhic fiber bundles are considered
 (topologicaly trivial) as being connections over the trivial fiber bundle,
 the symplectic form can be considered as being given by the
 natural symplectic form of the space of unitary connections:
 $$ \omega(\alpha_1 , \alpha_2 ) = \int_{\Si} tr(\alpha_1 \wedge \alpha_2) .$$
 Where $ \alpha_1 $ and $\alpha_2 $ are some $1$-forms with values in the 
 anti-symmetric endomorphisms of the trivial fiber bundle
 considered as elements of the tangent space of the space of connections.
 The space of moduli of the fiber bundles ${\cal M}_n(\Si)$
 can be considered as a symplectic quotient of the space of connections
 and the tangent space is then identified with the harmonic $1$-forms
 with values in the antisymmetric endomorphisms of the fiber bundle:
 $$ H^1(\Si, A({\mathbb C}^n)).$$
\subsection{The reduction of a holomorphic $1$-cocycle in a $1$-form over the
 surface $\Si$}
 There is no actual canonical and explicite mean to associate
 a connection with a holomorphic fiber bundle (stable). The
 problem is equivalent to find a metric of the stable holomorphic fiber
 bundle which allows to show him as an unitary representation
 of the Poincar\'e group. But, at the level of the tangent spaces
 in the trivial fiber bundle, the equivalent of the identification is simple,
 it is the isomorphism between the cohomology of \v Cech and the cohomology of 
 de Rham obtained via the Hodge theory \cite{Ho}.
 The reduction in a harmonic $1$-form is then:
 let a holomorphic $1$-cocycle $f$ be at the point $p$ and a choice of a 
partition of the unity adapted with the open sets
 of the surface which are the disc $D$ and the 
 surface $\Si$ minus the point $p$. There is then a  function $\varrho$
 with value $1$ near the point and zero outside.
 Then the form of type $(0,1)$ is considered:
 $ f \overline \partial \varrho$
  reduced  in its harmonic counterpart in the Dolbeault-Grothendieck 
 cohomology \cite{La}:
 $$ \phi (Q) = \int_{P \in \Si} f(P) \overline \partial_P \varrho(P) \wedge
 \partial_P \overline \partial_Q g(P,Q) .$$
 Where $g$ is the Green function of the Riemann surface $\Si$. The harmonic
 form of the de Rahm cohomology is then: $ \alpha = Re( \phi).$
\subsection{The reduction of the symplectic form of the double quotient}
 Let two holomorphic fiber bundles be $f_1$ and $f_2$ considered as 
 meromorhic functions over the disc $D$ and a pole in zero.

\medskip
\noindent
{\defi .
 $:g(P,Q):$ is the {\it renormalisation} of the Green function de Green:
 $$ :g(P,Q): = \  g(P,Q) - ln( |P-Q|),$$
 where $P , Q \in \Si$ and $|P-Q|$ is the geodesic distance over 
 the Riemann surface $\Si$.}

\medskip
\noindent
{\lem .  
 The Green function being biharmonic, it can be written
 $\partial_z \overline \partial_t :g(z,t):$ 
 as double series holomorphic, anti-holomorphic.}

\medskip
\noindent
{\it Proof}:
 the Green function de Green being biharmonic is the solution of the Laplacian
  $\Delta = \partial \overline \partial = \overline \partial \partial .$
 The partial derivations can then be written over the disc as
 holomorphic function in one of the variables and anti-holomorphic
 with the other; it shows, by the Cauchy formula, that the expression can be
 developped near zero in a double series.
 $$\partial_z \overline \partial_{t} :g(z,t):=
 \sum_{ n= k}^{+\infty} \sum_{m= p}^{+\infty}
 a_{n,m} z^n \overline t^m .$$
\medskip
 The symplectic form is then given by:

\medskip
\noindent
{\th .   \label{symp}
 The symplectic form corresponding with the tangent space of the 
 double quotient at level of the trivial fiber bundle is:
\begin{equation}
 \omega_e (f^1,f^2)=  (2 \pi)^2 Re ( \sum_{n}
 \sum_{m} a_{n,m}  tr( [f^1_{n-1}]^{*} f^2_{m-1}) ),
\end{equation}
 for $f_1 (z)= \sum_{r=q}^{+\infty} f^1_r z^r$ and 
 $f_2 (t) = \sum_{l=s}^{+\infty} f^2_l t^l$.}

\medskip
\noindent
{\it Proof}:
 the symplectic form is the real part of:
$$\int_{Q \in \Si}  tr(( \int_{P \in \Si} f_1^* (P) \partial_P \varrho(P)
 \wedge [ \overline 
\partial_P \partial_Q g(P,Q) ) ] \wedge $$ 
$$ \wedge (\int_{P \in \Si} f_2(P)
 \overline \partial_P \varrho(P) \wedge [ \partial_P
 \overline \partial_Q g(P,Q) ])). $$ 
There is a function $\varrho$ with $1$ value near the point and
 zero outside of the disc, the function goes towards
 the caracteristic function of the disc.
  Stockes is aplied two times and 
 the property of the function of
 Green is used \cite{La}:
$$ \int_{Q \in \Si}
 [ \partial_P \overline \partial_Q g(P,Q)] \wedge [ \partial_Q
 \overline \partial_{Q'} 
 g(Q,Q')] = \partial_P \overline \partial_{Q'} g(P,Q'). $$
 Instead an integration over all the surface $\Si$, 
 minus where $\varrho$ is $1$, as $\varrho$ goes towards the
 caracteristic function of the disc, Stockes is applied to $[\Si-D]$.
The following expression of the symplectic form is obtained:
$$ \int_{P \in S^1}  \int_{ Q \in S^1} tr(f_1^*(P) f_2(Q))
 [ \partial_P \overline \partial_Q : g(P,Q) :] .$$
 $\partial_P \overline \partial_Q :g(P,Q):$ is supposed being
 developped in a double series over the disc or a smaller;
 this double integrale is explained around zero:
 $$ \int \int_{z  , t \in S^1}  tr(f_1^*(z) f_2(t)) 
 \sum_{ n= k}^{+\infty} \sum_{m= p}^{+\infty}
 a_{n,m} z^n \overline t^m dz d\overline t.$$
 And so,
$$ \sum_{ n= k}^{+\infty}
 \sum_{m= p}^{+\infty} \sum_{r= q}^{+\infty} \sum_{l= s}^{+\infty}
 a_{n,m}  tr (  [f^1_r]^* f^2_l)  \int \int_{z  , t \in S^1}
 z^n \overline t^m \overline z^r t^l dz d\overline t ,$$
 developping in series the functions $f_1$ and $f_2$. So:
 $$  (2 \pi)^2\sum_{ n= k}^{+\infty}
 \sum_{m= p}^{+\infty} \sum_{r= q}^{+\infty} \sum_{l= s}^{+\infty}
 a_{n,m} tr([ f^1_r]^{*} f^2_l) \delta_{n-1,r} \delta_{m-1,l}, $$
 with $\delta$, the symbol of Kronecker.
\section{Conclusion}
 A presentation of the moduli spaces,
 ${\cal M}_n(\Si)$, over a Riemann surface $\Si$ has been done,
 and it has been showed how to obtain them in three different ways.
 A formula for the derivative of the regularised determinant has been computed.
 The existence of a simple formula has been showed 
 for the symplectic structure at level
 of the tangent space at the trivial fiber bundle of the double quotient.
\section{Acknowledgements}
I am greatly indebted to P.Gauduchon for his help, and kind advices.
\newpage
{ \scriptsize
\appendix
\section{The prodistributions}
\subsection{The product of two distributions}
 Let a compact differentiable variety  $M$ be, with orientation,
 the dual ${\cal D}'(M)$ \cite{Sc} of the  
 functions ${\cal C}^{{\infty}}$ is considered with compact support.

{\defi .
   ${\mathbb D}'_{2}(M)$, is the {\it space of $2$-distributions},
 the following topological tensor product \cite{G}:
$$ 
 {\mathbb D}'_{2} (M) :=  {\cal D}'(M) \hat \otimes_{{\cal C}^{{\infty}}(M)}
  {\cal D}'(M).
$$}
 This space, which is a  module over the ring of functions
   ${\cal C}^{{\infty}}$, allows a definition of a product of two
 distributions in the following sens:\\
 The product of two distributions $\eta$, $\eta'$ is the element
 $\eta \eta'  $:
$$
 \eta \eta' = \eta . \eta' = \eta \otimes \eta' \in {\mathbb D}'_{2}(M).
$$
\subsection{The product of $n$ distributions}
{\defi .
  ${\mathbb D}'_{n}(M)$ is the {\it space of $n$-distributions},
 the following tensor product:

$$ 
 {\mathbb D}'_{n} (M) :=  {\cal D}'(M) \hat \otimes_{{\cal C}^{{\infty}}(M)}
  {\cal D}'(M) \hat \otimes_{{\cal C}^{{\infty}}(M)}... \hat 
\otimes_{{\cal C}^{{\infty}}(M)} {\cal D}'(M) .
$$}
e over the ring of functions
   ${\cal C}^{{\infty}}$, allows a definition of a product of two
 distributions in the following sens:\\
 The product of two distributions $\eta$, $\eta'$ is the element
 $\eta \eta'  $:
$$
 \eta \eta' = \eta . \eta' = \eta \otimes \eta' \in {\mathbb D}'_{2}(M).
$$
\subsection{The product of $n$ distributions}
{\defi .
  ${\mathbb D}'_{n}(M)$ is the {\it space of $n$-distributions},
 the following tensor product:

$$ 
 {\mathbb D}'_{n} (M) :=  {\cal D}'(M) \hat \otimes_{{\cal C}^{{\infty}}(M)}
  {\cal D}'(M) \hat \otimes_{{\cal C}^{{\infty}}(M)}... \hat 
\otimes_{{\cal C}^{{\infty}}(M)} {\cal D}'(M) .
$$}

 The  tensor product is taken $n$ times.
 This space, a module over the ring of functions
  ${\cal C}^{{\infty}}$, allows 
 a definition of the product of $n$ distributions in the following sens:\\
 The product of $n$ distributions  $\eta_{1}$, $\eta_{2}$, ...,
$\eta_{n}$ is the element $\eta_{1}\eta_{2} ... \eta_{n} $:
$$
 \eta_{1}\eta_{2} ... \eta_{n}= \eta_{1} . \eta_{2} .\; ... \; .\eta_{n}
 = \eta_{1} \otimes \eta_{2} \otimes  ...  \otimes \eta_{n} \in
  {\mathbb D}'_{n}(M).
$$
\subsection{The space ${\mathbb D}'(M)$ of prodistributions}
$$
 {\mathbb D}'_{n}(M). {\mathbb D}'_{m}(M) \rightarrow  {\mathbb D}'_{n+m}(M).
$$
{\defi .
  ${\mathbb D}'(M)$, module
 over ${\cal C}^{{\infty}}(M)$, {\it space of prodistributions} is the
 following one:
$$
 {\mathbb D}'(M):=  lim_{n \in \NN, n \geq 2}  \; {\mathbb D}'_{n}(M).
$$}
 The limit is taken with the system of the inclusions \cite{G}.

\end{document}